\documentclass[11pt]{amsart}

\usepackage[T1]{fontenc}
\usepackage[utf8]{inputenc}
\usepackage{lmodern}
\usepackage{amsmath,amssymb,amsthm,mathtools}
\usepackage{graphicx}
\usepackage{float}
\usepackage{microtype}
\usepackage[colorlinks=true,linkcolor=blue,citecolor=blue,urlcolor=blue]{hyperref}

\newtheorem{theorem}{Theorem}[section]
\newtheorem{lemma}[theorem]{Lemma}

\theoremstyle{definition}
\newtheorem{definition}[theorem]{Definition}

\theoremstyle{remark}

\title[Support geometry of a parabolic nonlinear AR(1) process]{Support geometry of the stationary distribution for a parabolic nonlinear AR(1) process with binary innovations}
\author{Mukhriddin Khomidov}
\address{V.I. Romanovskiy Institute of Mathematics of the Uzbekistan Academy of Sciences, University Street 4, Olmazor District, Tashkent, Uzbekistan; International Islamic Academy of Uzbekistan, Abdulla Qodiriy Street 78, Tashkent, Uzbekistan}
\email{m.homidov1008@gmail.com}
\date{}

\begin{document}
\begin{abstract}
This paper studies the nonlinear AR(1) process
\(X_{n+1}=f(X_n)+\sigma\xi_{n+1}\), \(n\ge0\), driven by binary innovations
\(\xi_n\in\{0,1\}\), where \(f\) is strictly increasing and strictly concave,
with \(f(0)=0\), \(f'(0)=1\), and \(f'(x)<1\) for \(x>0\). The process is
shown to have a unique stationary distribution. The geometry of its support is
described in terms of the noise amplitude \(\sigma\) and the function \(f\):
the support is either an interval or a Cantor set of Lebesgue measure zero,
and the two regimes are separated by the relation between \(\sigma\) and the
unique number \(c_\sigma>0\) solving
\(c_\sigma-f(c_\sigma)=\sigma\). Namely, the support is an interval when
\(c_\sigma\ge2\sigma\) and a Cantor set when \(c_\sigma<2\sigma\). Under the
additional condition \(\ell_\infty:=\lim_{x\to\infty}f'(x)<1/2\), there
exists a unique critical value \(\sigma_c>0\), characterized by
\(f(2\sigma_c)=\sigma_c\), such that the support is an interval for
\(0<\sigma\le\sigma_c\) and a Cantor set for \(\sigma>\sigma_c\); if
\(\ell_\infty\ge1/2\), the support is an interval for every \(\sigma>0\).
In the Cantor regime the stationary distribution is non-atomic, and hence
singular continuous.
\end{abstract}
\maketitle

\section{Introduction}
Let \((\Omega, \mathcal{F}, P)\) be a probability space, and let
\(\{\xi_n\}_{n\geq 1}\) be a sequence of independent and identically
distributed random variables defined on it. Assume that \(X_0\geq0\)
is a random variable defined on this space and
independent of the sequence \(\{\xi_n\}_{n\geq 1}\).
Let \(f:\mathbb{R}_{+}\to\mathbb{R}_{+}\) be a nonnegative-valued function, where
\(\mathbb{R}_{+}=[0,+\infty)\). Introduce the following stochastic process:
\begin{equation}\label{eq1}
X_{n+1}=f(X_n)+\sigma\,\xi_{n+1},\qquad n\geq 0.
\end{equation}
Here \(\sigma>0\).
The stochastic process defined by relation (\ref{eq1}) is called a first-order
autoregressive process and is denoted by AR(1). Usually, the case where
\(f\) is linear is called linear autoregression, while the case where
\(f\) is nonlinear is called nonlinear autoregression.
Classical linear autoregressive processes are among the fundamental topics of time
series analysis and find wide application in practical problems. For the classical linear AR(1) recursion
$$
X_{n+1}=aX_n+\varepsilon_{n+1},\qquad |a|<1,
$$
the stationary solution admits the infinite moving-average representation
$$X_n=\sum_{k=0}^{\infty}a^k\varepsilon_{n-k},$$ so that its distribution is
completely determined by $a$ and the law of the innovations. When the
innovations possess smooth densities, the stationary law is typically studied
through its density and the corresponding integral equation. The situation
changes substantially for discrete innovations. Even in the linear AR(1)
setting, binary innovations may lead to Bernoulli-convolution-type stationary
laws, for which singularity and fractal support can occur. 

Nonlinear autoregressive processes have been studied in the works of many
researchers (e.g., \cite{Andel1989}, \cite{AndelHrach2000},
\cite{AndelGomezVega1989}, \cite{AndelNetukaZvara1984},
\cite{DieboltGuegan1993}, \cite{Jones1978}, and others).
Jones \cite{Jones1978} studied recursions of the form
\[
X_{n+1}=f(X_n)+\xi_{n+1},
\]
where \(f\) is a nonlinear regression function, and developed methods for
computing stationary distributions, moments, and predictive distributions,
along with numerical approximation techniques.
 Threshold and related nonlinear autoregressions were further studied by And\v{e}l, Netuka and Zv\'ara \cite{AndelNetukaZvara1984}, while And\v{e}l
\cite{Andel1989} gave a broader treatment of nonlinear time-series models and
their stationary distributions. Especially close to the present setting is
the work of And\v{e}l, G\'omez and Vega \cite{AndelGomezVega1989}, where nonlinear AR(1) processes with two-point
innovations were studied for a particular power-law regression
\[
f(x)=x^{1/h},\qquad \quad h>1,
\]
together with procedures for computing stationary distributions and moments.
Later, And\v{e}l and Hr\'ach \cite{AndelHrach2000} developed iterative
methods for calculating stationary densities, while Diebolt and Gu\'egan
\cite{DieboltGuegan1993} investigated tail properties of stationary
densities for general nonlinear autoregressions of order one.
And\v{e}l, G\'omez and Vega \cite{AndelGomezVega1989} proved the existence of
a stationary distribution for a nonlinear autoregressive process generated by
a power function and driven by binary innovations. However, they did not
establish whether this stationary distribution is absolutely continuous,
singular, or of mixed type. A natural question is thus whether the stationary
distribution belongs to the absolutely continuous, singular, or mixed type.

One of the main motivations of this work is the question of which type the
stationary distribution belongs to.

Suppose the function \(f: \mathbb{R}_+\to \mathbb{R}_+\) satisfies the following conditions:
\begin{enumerate}
    \item \(f\in C^{2}([0,\infty))\);
    \item \(f(0)=0\) and \(f'(0)=1\);
    \item \(f'(x)>0\) for all \(x>0\);
    \item \(f''(x)<0\) for all \(x\geq 0\).
\end{enumerate}
The above conditions imply that \(f\) is concave and strictly increasing, and
that \(f'(x)<1\) for all \(x>0\). Consequently, \(f\) is a contraction on
every interval \([a,b]\subset(0,\infty)\).

In this paper, we consider process (1) under the assumption that the function
\(f\) satisfies the above conditions. Assume that \(\{\xi_n\}_{n\geq1}\) is an i.i.d.\ Bernoulli sequence satisfying
\[
\mathbb P(\xi_n=0)=p,\qquad
\mathbb P(\xi_n=1)=q:=1-p,\qquad
0<p<1.
\]
The main purpose of this paper is to study the stationary distribution of the
nonlinear process (1) generated by the above function \(f\) and driven by
Bernoulli innovations, together with the support of this distribution and its
geometry.

\section{Notation}
Let \((\mathcal{X},\rho)\) be a metric space. For a continuous function
\(h:\mathcal{X}\to\mathbb{R}\), denote by \(\operatorname{Lip}(h)\) the
Lipschitz constant defined as
\[
\operatorname{Lip}(h):=\sup_{x\neq y}\frac{|h(x)-h(y)|}{\rho(x,y)}.
\]
Let \(\mathcal P_1(\mathcal X)\) denote the probability measures on
\(\mathcal X\) with finite first moment. For
\(\mu_1,\mu_2\in\mathcal P_1(\mathcal X)\), the 1-Wasserstein distance is
\[
W_1(\mu_1,\mu_2)
=\inf_{\substack{Y\sim\mu_1\\ Z\sim\mu_2}}
\mathbb{E}\bigl[\rho(Y,Z)\bigr],
\]
where the infimum is taken over all couplings of \(\mu_1\) and \(\mu_2\).
On the present state space \([0,\infty)\), the Kantorovich--Rubinstein dual
formula gives
\[
W_1(\mu_1,\mu_2)
=\sup_{\operatorname{Lip}(h)\le1}
\left|\int h\,d\mu_1-\int h\,d\mu_2\right|;
\]
see, for example, \cite{Villani2009}.
Let \(T_i:\mathcal{X}\to\mathcal{X}\), \(i=1,\ldots,m\), be functions, and let
\(\Phi=\{T_1,\ldots,T_m\}\) be an iterated function system. Consider the
probability space \((\Xi,\mathcal{A},P)\), where
\[
\Xi=\Bigl\{\underline{\omega}=(\omega_1,\omega_2,\ldots,\omega_n,\ldots):
\ \omega_i\in\{0,1,\ldots,m-1\}\Bigr\},
\]
\(\mathcal{A}\) is the \(\sigma\)-algebra generated by the cylinder sets in
\(\Xi\). We also define
\[
\Xi_n=\{0,1,\ldots,m-1\}^{n}.
\]
Denote by
\[
F_n(x)=T_{\omega_1}\circ T_{\omega_2}\circ\cdots\circ T_{\omega_n}(x)
\]
and
\[
\widetilde{F}_n(x)
=T_{\omega_n}\circ T_{\omega_{n-1}}\circ\cdots\circ T_{\omega_1}(x)
\]
the two different compositions of these mappings. It is clear that, if the coordinate functions\(\eta_n:=\eta_n(\underline{\omega})=\omega_n\)
are independent and identically distributed, then the random variables
\(F_n(X_0)\) and \(\widetilde{F}_n(X_0)\) have the same distribution. 

Let \(X_0\) be an \(\mathcal{X}\)-valued random variable. Then
\(
X_n=\widetilde{F}_n(X_0)
\)
is a Markov chain. (see. \cite{David}).

The following theorem due to Letac plays an important role in the study of
stationary distributions of iterated function systems.
\begin{theorem}[\cite{Letac}]
Assume that the distribution \(\nu\) is concentrated on continuous functions
and that \(F_\infty(x)=\lim\limits_{n\to\infty}F_n(x)\)
exists almost surely and is independent of \(x\). Then the distribution of
\(
X_\infty=F_\infty(x)
\) is the unique stationary distribution for the Markov chain
\(\widetilde{F}_n(x)\). Moreover, it is attractive: any compactly supported
initial distribution converges to it under the action of the Markov operator.
\end{theorem}

If \(\operatorname{Lip}(T_i)<1\) for all \(i=0,\ldots,m-1\), then the system is called a contractive iterated function system. Contractive systems of this type were thoroughly studied by Hutchinson \cite{Hutchinson1981}. Hutchinson's classical theory \cite{Hutchinson1981} treats uniformly contractive iterated function systems. Indeed, if every mapping \(T_i\), \(i=0,\ldots,m-1\), were Lipschitz continuous with \(\operatorname{Lip}(T_i)<1\), then each \(T_i\) would be a strict contraction, and the basic existence and uniqueness results for the invariant set and stationary probability measure would follow from the classical contractive theory.

The existence and uniqueness of the stationary distribution subsequently led
to the study of its support and regularity. Bagaev and Ganeeva
\cite{BagaevGaneeva2025} obtained conditions ensuring Cantor attractors for
uniformly contractive iterated function systems. However, Cantor topology
alone does not imply singularity with respect to Lebesgue measure, since
Cantor sets of positive Lebesgue measure exist.

\section{Main results}
The main results of the paper are stated next.
\begin{theorem}\label{thm3.1}
For every $\sigma>0$, the process \eqref{eq1} admits a unique
stationary probability measure $\mu_\sigma$.  Moreover, for every
initial distribution of $X_0$ on $[0,\infty)$,
\(
\mathcal L(X_n)\Longrightarrow\mu_\sigma
\qquad\text{as }n\to\infty,
\)
where $\mathcal L(X_n)$ denotes the law of $X_n$ and
$\Longrightarrow$ denotes weak convergence of probability measures.
\end{theorem}
\begin{proof} Set \(\mathcal{X}=\mathbb{R}_{+}\) and define the two maps
\[
T_0(x):=f(x),
\qquad
T_1(x):=f(x)+\sigma.
\]
Then the recursion \eqref{eq1} can be written as
\[
X_{n+1}=T_{\xi_{n+1}}(X_n),
\qquad n\geq0.
\]
Consequently,
\[
X_n
=
T_{\xi_n}\circ T_{\xi_{n-1}}\circ\cdots\circ T_{\xi_1}(X_0).
\]

The corresponding Markov operator acting on probability measures is
\[
\mathcal P\mu
=
p\,(T_0)_*\mu
+
q\,(T_1)_*\mu,
\]
where \(p=\mathbb P(\xi_1=0)\) and \(q=\mathbb P(\xi_1=1)=1-p\).
A probability measure \(\mu\) is called stationary if
\(\mathcal P\mu=\mu.\)

Since \(f''(x)<0\) for \(x\geq0\), the derivative \(f'\) is strictly decreasing. Together with \(f'(0)=1\) and \(f'(x)>0\) for \(x>0\), this gives
\[
0<f'(x)\leq1,
\qquad x\geq0.
\]
In particular,
\[
r_\sigma:=f'(\sigma)<1.
\]
Moreover, since
\[
T_0'(x)=T_1'(x)=f'(x),
\]
both maps are non-expanding on \([0,\infty)\).
Let
\[
r_\sigma:=f'(\sigma)<1,\qquad
\rho_\sigma:=p+q\,r_\sigma<1,
\]
where \(p=P(\xi_1=0)\) and \(q=P(\xi_1=1)\), with \(p+q=1\) and \(q>0\).

Let \(\widetilde{F}_n(x)\) and \(\widetilde{F}_n(y)\) be two copies of the
process starting from \(x,y\ge 0\), respectively, and driven by the same
Bernoulli sequence. 

For the first step we have
\[
|\widetilde{F}_1(x)-\widetilde{F}_1(y)|
=|T_{\xi_1}(x)-T_{\xi_1}(y)|
\le f'(0)\,|x-y|=|x-y|.
\]
For the second step, consider the two possible values of \(\xi_1\):

- If \(\xi_1=0\), then
  \[
  |\widetilde{F}_2(x)-\widetilde{F}_2(y)|
  =|T_{\xi_2}(f(x))-T_{\xi_2}(f(y))|
  \le f'(0)\,|f(x)-f(y)|\le |x-y|.
  \]
- If \(\xi_1=1\), then
  \(\widetilde{F}_1(x)=f(x)+\sigma\ge\sigma\) and
  \(\widetilde{F}_1(y)=f(y)+\sigma\ge\sigma\). Since \(f'\) is decreasing,
  \[
  |\widetilde{F}_2(x)-\widetilde{F}_2(y)|
  \le r_\sigma\,|\widetilde{F}_1(x)-\widetilde{F}_1(y)|
  \le r_\sigma\,|x-y|.
  \]

Taking expectation over \(\xi_1\),
\[
\mathbb{E}\bigl|\widetilde{F}_2(x)-\widetilde{F}_2(y)\bigr|
\le \rho_\sigma\,|x-y|,
\qquad x,y\ge 0.
\]

Let
\[
P_n=\mathcal{L}\bigl(\widetilde{F}_n(X_0)\bigr).
\]

We now show that \(\mathcal{P}^2\) is a contraction on
\(\mathcal{P}_1([0,\infty))\) with respect to \(W_1\). Let \(\mu,\nu\in
\mathcal{P}_1([0,\infty))\), and let \((X,Y)\) be an optimal coupling of
\((\mu,\nu)\), so that \(\mathbb{E}|X-Y|=W_1(\mu,\nu)\). The two-step
estimate above yields
\[
\mathbb{E}\Bigl[
\bigl|\widetilde{F}_2(X)-\widetilde{F}_2(Y)\bigr|
\,\Big|\,X,Y\Bigr]
\le \rho_\sigma\,|X-Y|.
\]
Taking expectations and using the definition of \(W_1\),
\[
W_1(\mathcal{P}^2\mu,\mathcal{P}^2\nu)
\le \mathbb{E}\bigl|\widetilde{F}_2(X)-\widetilde{F}_2(Y)\bigr|
\le \rho_\sigma\,W_1(\mu,\nu).
\]
Hence \(\mathcal{P}^2\) is a contraction on the complete metric space
\(\bigl(\mathcal{P}_1([0,\infty)),W_1\bigr)\).

Next we verify that \(\mathcal{P}\) maps \(\mathcal{P}_1([0,\infty))\)
into itself. Since \(f\) is concave, \(f(0)=0\), and \(f'(0)=1\), we have
\(f(x)\le x\) for all \(x\ge 0\), and therefore
\[
\int x\,d(\mathcal{P}\mu)
=p\int f(x)\,d\mu+q\int\bigl(f(x)+\sigma\bigr)\,d\mu
\le \int x\,d\mu+q\sigma<\infty.
\]

By the Banach fixed-point theorem, \(\mathcal{P}^2\) has a unique fixed
point \(\mu_\sigma\in\mathcal{P}_1([0,\infty))\), i.e.\
\(\mathcal{P}^2\mu_\sigma=\mu_\sigma\). Since
\[
\mathcal{P}^2(\mathcal{P}\mu_\sigma)
=\mathcal{P}(\mathcal{P}^2\mu_\sigma)
=\mathcal{P}\mu_\sigma,
\]
the uniqueness of the fixed point of \(\mathcal{P}^2\) implies
\(\mathcal{P}\mu_\sigma=\mu_\sigma\). Thus \(\mu_\sigma\) is a stationary
distribution of the process.

If the initial law \(\nu\) belongs to
\(\mathcal P_1([0,\infty))\), then the two-step contraction and the
one-step non-expansiveness give the quantitative estimate
\[
W_1(\mathcal P^n\nu,\mu_\sigma)
\le \rho_\sigma^{\lfloor n/2\rfloor}W_1(\nu,\mu_\sigma),
\]
and hence \(\mathcal P^n\nu\Rightarrow\mu_\sigma\).

It remains to justify the asserted weak convergence for an arbitrary initial
probability law \(\nu\), which need not have a finite first moment. For every
fixed \(x\ge0\), the preceding estimate applied to \(\delta_x\) gives
\[
\mathcal P^n\delta_x\Rightarrow\mu_\sigma.
\]
Let \(\varphi\) be bounded and continuous. Then
\[
\int \varphi\,d(\mathcal P^n\nu)
=\int_{[0,\infty)}
\left(\int \varphi\,d(\mathcal P^n\delta_x)\right)\nu(dx).
\]
For each \(x\), the inner integral converges to
\(\int\varphi\,d\mu_\sigma\), while its absolute value is bounded by
\(\|\varphi\|_\infty\). The dominated convergence theorem therefore yields
\[
\mathcal P^n\nu\Rightarrow\mu_\sigma.
\]

Finally, if \(\nu\) is any stationary probability measure, then
\(\mathcal P^n\nu=\nu\) for every \(n\). Applying the weak convergence just
proved to the initial law \(\nu\) gives \(\nu=\mu_\sigma\). Thus the
stationary probability measure is unique. This completes the proof of
Theorem~\ref{thm3.1}.
\end{proof}

Let \(K_\sigma\) be the support of the stationary distribution, that is,
\(
K_\sigma=\operatorname{supp} \mu_\sigma .
\)
The geometry of the stationary support is governed by the unique number
\(c_\sigma>0\) defined by
\begin{equation}\label{eq2}
c_\sigma-f(c_\sigma)=\sigma.
\end{equation} 
Now we formulate following useful lemmas.
\begin{lemma}\label{lem3.1}
For every $\sigma>0$, there exists a unique number $c_\sigma>0$
satisfying \eqref{eq2}. \end{lemma}

\begin{proof}
Set
$
h(x):=x-f(x),\qquad x\geq0.
$
Since $f(0)=0$, we have $h(0)=0$.  Moreover, strict concavity and
$f'(0)=1$ imply
$
0<f'(x)<1,\qquad x>0,
$
and hence
$
h'(x)=1-f'(x)>0.
$
Thus $h$ is strictly increasing.

Fix $a>0$.  Since $f'$ is decreasing and $f'(a)<1$, for $x\geq a$,
\[
h'(x)\geq1-f'(a)>0,
\]
so that $h(x)\to\infty$ as $x\to\infty$. 
Therefore the equation $h(x)=\sigma$, which is
equivalent to \eqref{eq2}, has a unique solution $c_\sigma>0$.
Finally,
$
c_\sigma=\sigma+f(c_\sigma)>\sigma.
$
\end{proof}
Define
$I_\sigma:=[0,c_\sigma].$
\begin{lemma}\label{lem3.3}
The interval $I_\sigma$ is invariant under both $T_0$ and $T_1$.
More precisely,
$
T_0(I_\sigma)
=
[0,c_\sigma-\sigma],
$
and
$
T_1(I_\sigma)
=
[\sigma,c_\sigma].
$
\end{lemma}

\begin{proof}
Since $f$ is strictly increasing and
$
f(c_\sigma)=c_\sigma-\sigma,
$
we have
\[
T_0(I_\sigma)
=
[f(0),f(c_\sigma)]
=
[0,c_\sigma-\sigma].
\]
Similarly,
\[
T_1(I_\sigma)
=
[\sigma,f(c_\sigma)+\sigma]
=
[\sigma,c_\sigma].
\]
Both intervals are contained in $I_\sigma$.
\end{proof}

For each \(\underline{\omega}\in\Xi\), consider the nested sets
\[
J_{\omega_1\cdots\omega_n}
:=T_{\omega_1}\circ\cdots\circ T_{\omega_n}(I_\sigma),
\qquad n\ge1.
\]
It is clear that
\[
J_{\omega_1}
\supset
J_{\omega_1\omega_2}
\supset
J_{\omega_1\omega_2\omega_3}
\supset\cdots .
\]
We need following lemma.

\begin{lemma}\label{lem3.4}
Let $ D_n:=\max_{w\in\Xi_n}\operatorname{diam}J_w$  Then, $ D_n\longrightarrow0
\qquad\text{as }n\to\infty. $
\end{lemma}\begin{proof}
For \(x\ge y\ge0\), the monotonicity of \(f'\) gives
\[
f(x)-f(y)
=\int_y^x f'(t)\,dt
=\int_0^{x-y}f'(y+s)\,ds
\le \int_0^{x-y}f'(s)\,ds
=f(x-y).
\]
Consequently, for \(i=1,2\),
\[
|T_i(x)-T_i(y)|=|f(x)-f(y)|\le f(|x-y|).
\]
If \(J\subset I_\sigma\) is an interval, it follows that
\[
\operatorname{diam}T_i(J)\le f(\operatorname{diam}J).
\]
Iterating this estimate along any word \(w\in\Xi_n\) yields
\[
\operatorname{diam}J_w\le f^n(c_\sigma).
\]
Since \(0<f(x)<x\) for every \(x>0\), the sequence
\(f^n(c_\sigma)\) decreases to a limit \(L\ge0\). By continuity,
\(L=f(L)\), and the strict inequality \(f(x)<x\) for \(x>0\) forces
\(L=0\). Hence
\[
D_n\le f^n(c_\sigma)\longrightarrow0,
\]
uniformly over all words of length \(n\).
\end{proof}
\begin{definition}\label{def3.4}
For
$
\underline{\omega}=(\omega_1,\omega_2,\ldots)\in\Xi,
$
define $\pi_\sigma(\underline{\omega})$ by
\[
\bigcap_{n=1}^{\infty}
J_{\omega_1\cdots\omega_n}
=
\{\pi_\sigma(\underline{\omega})\}.
\]
The map
$
\pi_\sigma:\Xi\longrightarrow I_\sigma
$
is called the coding map.
\end{definition}
By Lemma~\ref{lem3.4}, for each \(\underline{\omega}\) the set
\(\{\pi_\sigma(\underline{\omega})\}\) is a singleton. The same uniform
diameter estimate also shows that \(\pi_\sigma\) is continuous: if two
sequences agree in their first \(n\) symbols, their images lie in the same
cylinder interval of diameter at most \(D_n\).

Let \(S:\Xi\to\Xi\) be the left shift and let \(\beta_{p,q}\) be the
Bernoulli product measure on \(\Xi=\{0,1\}^{\mathbb N}\) with weights
\(p\) and \(q\), that is,
\(\beta_{p,q}(\omega_1=0)=p\) and
\(\beta_{p,q}(\omega_1=1)=q\). The coding relation
\[
\pi_\sigma(\underline\omega)
=T_{\omega_1}\bigl(\pi_\sigma(S\underline\omega)\bigr)
\]
follows directly from the definition. Hence the probability measure
\[
\nu_\sigma:=(\pi_\sigma)_*\beta_{p,q}
\]
satisfies, for every bounded measurable \(\varphi\),
\[
\int\varphi\,d\nu_\sigma
=p\int\varphi\circ T_1\,d\nu_\sigma
+q\int\varphi\circ T_2\,d\nu_\sigma.
\]
Thus \(\mathcal P\nu_\sigma=\nu_\sigma\), and the uniqueness in
Theorem~\ref{thm3.1} gives
\[
\mu_\sigma=(\pi_\sigma)_*\beta_{p,q}.
\]

Set \(K^*_\sigma:=\pi_\sigma(\Xi)\). Since \(\Xi\) is compact and
\(\pi_\sigma\) is continuous, \(K^*_\sigma\) is compact, and
\(\mu_\sigma(K^*_\sigma)=1\); hence
\(\operatorname{supp}\mu_\sigma\subseteq K^*_\sigma\). Conversely, if
\(x=\pi_\sigma(\underline\omega)\in K^*_\sigma\) and \(U\) is an open
neighbourhood of \(x\), continuity of \(\pi_\sigma\) gives a cylinder
neighbourhood \(C\) of \(\underline\omega\) such that
\(\pi_\sigma(C)\subset U\). Since \(p,q>0\),
\(\beta_{p,q}(C)>0\), and therefore \(\mu_\sigma(U)>0\). Thus
\(K^*_\sigma\subseteq\operatorname{supp}\mu_\sigma\). Consequently,
\[
\pi_\sigma(\Xi)=\operatorname{supp}\mu_\sigma=K_\sigma.
\]
\begin{theorem}
\label{thm3.6}
For every $\sigma>0$,
\[
K_\sigma=
\begin{cases}
[0,c_\sigma],
& c_\sigma\geq2\sigma,\\[4pt]
\text{a Cantor set},
& c_\sigma<2\sigma.
\end{cases}
\]
\end{theorem}
\begin{proof}
By the definition of $c_\sigma$,
$
f(c_\sigma)=c_\sigma-\sigma,
$
and therefore
\[
T_0(I_\sigma)=[0,c_\sigma-\sigma],
\qquad
T_1(I_\sigma)=[\sigma,c_\sigma].
\]
Assume first that
$
c_\sigma\geq 2\sigma.
$
Then
$
c_\sigma-\sigma\geq\sigma,
$
and hence
$
T_0(I_\sigma)\cup T_1(I_\sigma)=I_\sigma.
$
Let \(x\in I_\sigma\) be arbitrary. Since
\(I_\sigma=T_0(I_\sigma)\cup T_1(I_\sigma)\), there exist
\(\omega_1\in\{0,1\}\) and \(x_1\in I_\sigma\) such that
\(x=T_{\omega_1}(x_1)\). Iterating this argument yields a sequence
\(\underline{\omega}=(\omega_1,\omega_2,\ldots)\in\Xi\) for which
\[
x\in T_{\omega_1}\circ T_{\omega_2}\circ\cdots\circ T_{\omega_n}(I_\sigma)
=:J_{\omega_1\cdots\omega_n}
\qquad\text{for all }n\ge 1.
\]
By Lemma~\ref{lem3.4},
\(
\operatorname{diam}J_{\omega_1\cdots\omega_n}\longrightarrow 0
\qquad (n\to\infty),
\)
and hence
\[
x=\bigcap_{n\ge1} J_{\omega_1\cdots\omega_n}
=\pi_\sigma(\underline{\omega}).
\]
Thus $I_\sigma\subset K_\sigma$.  Since the reverse inclusion is
immediate from the construction,
\[
K_\sigma=I_\sigma=[0,c_\sigma].
\]

Suppose now that
$
c_\sigma<2\sigma.
$
Then
\[
T_0(I_\sigma)\cap T_1(I_\sigma)=\varnothing.
\]
Let $\underline{\omega},\underline{\eta}\in\Xi$ be distinct, and let $n$ be the first
index for which $\omega_n\neq\eta_n$.  Writing
\[
u=(\omega_1,\ldots,\omega_{n-1})
  =(\eta_1,\ldots,\eta_{n-1}),
\]
we have
\[
\pi_\sigma(\underline{\omega})\in T_u(T_{\omega_n}(I_\sigma)),
\qquad
\pi_\sigma(\underline{\eta})\in T_u(T_{\eta_n}(I_\sigma)).
\]
Since $T_u$ is injective and
$T_0(I_\sigma)\cap T_1(I_\sigma)=\varnothing$, these two sets are
disjoint.  Hence
\[
\pi_\sigma(\underline{\omega})\neq\pi_\sigma(\underline{\eta}),
\]
so that the coding map $\pi_\sigma$ is injective.

Since \(\Xi\) is a Cantor space and \(\pi_\sigma\) is a homeomorphism of
\(\Xi\) onto \(K_\sigma\), it follows that \(K_\sigma\) is a Cantor set.
\end{proof}

According to Theorem~\ref{thm3.6}, the support is either a Cantor set or an
interval, depending on the relationship between \(\sigma\) and
\(c_\sigma\). The following theorem examines this relationship for a given
function \(f\). Set
\(
\ell_\infty:=\lim_{x\to\infty}f'(x).
\)
\begin{theorem} \label{thm3.7}
If $\ell_\infty<1/2$, then there exists a unique critical value
$\sigma_c>0$ characterized by $ f(2\sigma_c)=\sigma_c.$
Moreover, $K_\sigma=[0,c_\sigma]
\qquad\text{ for }\qquad 0<\sigma\leq\sigma_c,
$
whereas \(K_\sigma \quad\text{is a Cantor set for }\sigma>\sigma_c.\)
If $\ell_\infty\geq1/2$, then \( K_\sigma=[0,c_\sigma] \qquad\text{for every }\sigma>0.\)
\end{theorem}
\begin{proof}
Let
$
h(x)=x-f(x).
$
As observed above,
\[
h'(x)=1-f'(x)>0,
\qquad
h''(x)=-f''(x)>0
\]
for $x>0$.  Thus $h$ is strictly increasing and strictly convex, and
$
c_\sigma=h^{-1}(\sigma).
$
Clearly, \(c_\sigma\) is a function of \(\sigma\). Since \(h\) is a strictly
increasing bijection of \([0,\infty)\) onto \([0,\infty)\), the relation
\[
c_\sigma-f(c_\sigma)=\sigma
\]
determines \(c_\sigma\) uniquely. Differentiating this identity with respect
to \(\sigma\), we get
\[
c_\sigma'=\frac{1}{1-f'(c_\sigma)}>0,
\qquad
c_\sigma''=\frac{f''(c_\sigma)}{\bigl(1-f'(c_\sigma)\bigr)^{3}}<0.
\]
Hence $\sigma\mapsto c_\sigma$ is strictly increasing and strictly
concave.  Since $c_0=0$, strict concavity implies that
$
\sigma\longmapsto\frac{c_\sigma}{\sigma}
$
is strictly decreasing on $(0,\infty)$.

Moreover, since
$
f'(0)=1,
$
we have\(\lim\limits_{x\to 0+}\frac{f(x)}{x}=1\)
and therefore
\[
\lim\limits_{x\to 0+}\frac{h(x)}{x}
=
1-\frac{f(x)}{x}=0.\]
Since $\sigma=h(c_\sigma)$, it follows that
$
\lim\limits_{\sigma\to 0+}\frac{c_\sigma}{\sigma}
=\lim\limits_{\sigma\to 0+}\frac{c_\sigma}{h(c_\sigma)}=\infty.
$
Since $f'$ is positive and decreasing, the limit
$
\ell_\infty=\lim\limits_{x\to\infty}f'(x)
$
exists.  Furthermore,
\[
f(x)=\int_0^x f'(t)\,dt,
\]
and hence
$
\lim\limits_{x\to\infty}\frac{f(x)}{x}=\ell_\infty.
$
Consequently,
$
\lim\limits_{\sigma\to\infty}\frac{\sigma}{c_\sigma}
=1-\ell_\infty,
$
and thus
$
\lim\limits_{\sigma\to\infty}\cfrac{c_\sigma}{\sigma}=
\frac{1}{1-\ell_\infty}.
$

Assume first that $\ell_\infty<1/2$.  Then
$
\frac{1}{1-\ell_\infty}<2.
$
Since $c_\sigma/\sigma$ is continuous and strictly decreasing from
$+\infty$ to a value strictly smaller than $2$, there exists a unique
$\sigma_c>0$ such that
\[
\frac{c_{\sigma_c}}{\sigma_c}=2,
\]
or equivalently,
$
c_{\sigma_c}=2\sigma_c.
$
Using
\[
c_{\sigma_c}-f(c_{\sigma_c})=\sigma_c,
\]
it follows that
\[
f(2\sigma_c)=\sigma_c.
\]

Moreover,
\[
c_\sigma>2\sigma
\quad\text{for }0<\sigma<\sigma_c,
\]
whereas
\[
c_\sigma<2\sigma
\quad\text{for }\sigma>\sigma_c.
\]
The conclusion now follows from
Theorem~\ref{thm3.6}.

If $\ell_\infty\geq1/2$, then
\[
\frac{1}{1-\ell_\infty}\geq2.
\]
Since $c_\sigma/\sigma$ is strictly decreasing to this limit, we have
\[
\frac{c_\sigma}{\sigma}>2
\qquad\text{for every finite }\sigma>0.
\]
Hence
$
c_\sigma>2\sigma
$
for every $\sigma>0$, and Theorem~\ref{thm3.6} yields
$
K_\sigma=[0,c_\sigma].
$
\end{proof}

When the support is of Cantor type, it is important to determine its
Lebesgue measure. If the Lebesgue measure of the support is zero, then the
stationary measure is singular with respect to Lebesgue measure. In what
follows, we analyze precisely the Cantor-type support. 

The following theorem is the main result of this paper.
\begin{theorem}  \label{thm3.8}
If $c_\sigma<2\sigma,$
then \(\lambda(K_\sigma)=0,\) and \(\mu_\sigma\) has no atoms. Hence
\(\mu_\sigma\perp\lambda\), and \(\mu_\sigma\) is singular continuous
with respect to Lebesgue measure. \end{theorem} 
Here $\lambda$ denotes Lebesgue measure. 
\begin{proof}
Assume that $
c_\sigma<2\sigma,
$
and write, for simplicity,
\( c:=c_\sigma, \qquad I=[0,c]. \)
Then
\[
T_0(I)=[0,c-\sigma],
\qquad
T_1(I)=[\sigma,c],
\]
and
$
c-\sigma<\sigma.
$
Hence the first-level cylinders are disjoint.  
By Theorem~\ref{thm3.6}, the coding map $\pi_\sigma$ is injective.
We first prove that
$
\lambda(K_\sigma)=0.
$

Let $u$ be a finite word and write
\(T_u=T_{u_1}\circ\cdots\circ T_{u_m}.\)
Each $T_i$ is increasing and concave, and therefore so is every
composition $T_u$.

Divide the interval \(J_u=T_u(I)\) into the following three pairwise
disjoint parts:
\[
J_{u0}=T_u([0,c-\sigma]),\qquad
J_{u1}=T_u([\sigma,c]),\qquad
G_u=\bigl(T_u(c-\sigma),T_u(\sigma)\bigr).
\]
Since \(T_u\) is increasing and \(c-\sigma<\sigma\), it is easy to see that
\[
J_u=J_{u0}\cup J_{u1}\cup G_u
\]
is a disjoint union. Moreover, \(G_u\) is a gap in \(K_\sigma\).
Each \(T_i\) is increasing and concave. Since the composition of a concave
increasing function with a concave function is again concave, every
composition \(T_u\) is concave. Hence 
\[
\frac{T_u(\sigma)-T_u(c-\sigma)}
     {2\sigma-c}
\geq
\frac{T_u(c)-T_u(\sigma)}
     {c-\sigma}.
\]
Therefore
\begin{equation}\label{eq3}
|G_u|
\geq
\kappa\,|J_{u1}|,
\qquad
\kappa:=
\frac{2\sigma-c}{c-\sigma}>0.
\end{equation}

Let $x\in K_\sigma$ have symbolic code
$
\underline{\omega}=(\omega_1,\omega_2,\ldots)
$
containing infinitely many symbols equal to $1$.  For every index $n$
such that $\omega_n=1$, let
\(u=(\omega_1,\ldots,\omega_{n-1}).
\)
Then $ x\in J_{u1}.$
Consider the one-sided interval
\[
H_u
:=
[T_u(c-\sigma),x].
\]
Since $G_u\subset H_u\setminus K_\sigma$ and
\[
|H_u|
\leq
|G_u|+|J_{u2}|,
\]
\eqref{eq3} yields
\[
\frac{\lambda(H_u\setminus K_\sigma)}{|H_u|}
\geq
\frac{\kappa}{1+\kappa}
=
\frac{2\sigma-c}{\sigma}
>0.
\]
Moreover,
$
H_u\subset J_u,
$
and the diameters of the cylinder intervals tend to zero.  Hence
$
|H_u|\longrightarrow0
$
along the subsequence of indices for which $\omega_n=2$.

Let $r_u:=|H_u|$. Since
$
H_u\subset [x-r_u,x+r_u],
$
we obtain
\[
\frac{\lambda([x-r_u,x+r_u]\setminus K_\sigma)}{2r_u}
\geq
\frac{1}{2}\,\frac{2\sigma-c}{\sigma}>0.
\]
As $r_u\to0$ along the above subsequence, such an $x$ cannot be a
Lebesgue density point of $K_\sigma$.
The only remaining points are those whose symbolic codes contain only
finitely many $1$'s.  Every such code has the form
\[
u000\cdots
\]
for some finite word $u$, and hence these points form at most a
countable set.

If $\lambda(K_\sigma)>0$, the Lebesgue density theorem would imply that
almost every point of $K_\sigma$ is a density point.  The preceding
argument shows that all density points would have to belong to a
countable set, which is impossible.  Consequently,
$\lambda(K_\sigma)=0.$
It is clear that,
$\mu_\sigma(K_\sigma)=1,
$
and therefore
$$\mu_\sigma\perp\lambda.$$

It remains to exclude atoms. Since \(\pi_\sigma\) is injective and
\(\pi_\sigma(\Xi)=K_\sigma\), each \(x\in K_\sigma\) corresponds to a
unique code \(\underline{\omega}=(\omega_1,\omega_2,\ldots)\in\Xi\) with
\(\pi_\sigma(\underline{\omega})=x\). By the representation
\(\mu_\sigma=(\pi_\sigma)_*\beta_{p,q}\), and because distinct
level-\(n\) cylinder intervals are disjoint in the separated regime, the
level-\(n\) cylinder containing \(x\) has mass
\[
\mu_\sigma\bigl(J_{\omega_1\cdots\omega_n}\bigr)
=\prod_{i=1}^{n}
\begin{cases}
p,&\omega_i=0,\\
q,&\omega_i=1,
\end{cases}
\le \max\{p,q\}^{n}.
\]
Since \(\{x\}\subset J_{\omega_1\cdots\omega_n}\) for every \(n\ge1\),
\[
\mu_\sigma(\{x\})
\le \max\{p,q\}^{n}.
\]
As \(\max\{p,q\}<1\), letting \(n\to\infty\) yields
\(\mu_\sigma(\{x\})=0\). Thus \(\mu_\sigma\) has no atoms. Being singular
with respect to Lebesgue measure and atomless, \(\mu_\sigma\) is singular
continuous.
\end{proof}

Consider the process
\[
X_{n+1}=f(X_n)+\sigma\,\xi_{n+1}
\]
with \(X_0=0\). The function
\[
f(x)=\frac{x}{1+x}
\]
satisfies the above conditions. For this function,
\[
c_\sigma=\frac{\sigma+\sqrt{\sigma^2+4\sigma}}{2},
\qquad
\sigma_c=\frac12,
\]
so that the critical value separating the two regimes is
\(\sigma_c=1/2\). In the case \(p=q=1/2\), for \(\sigma=0.2\) and
\(\sigma=1\), the plots of the empirical distribution function of
\(X_{100}\) are given below.

\begin{figure}[H]
\centering
\includegraphics[width=0.75\textwidth]{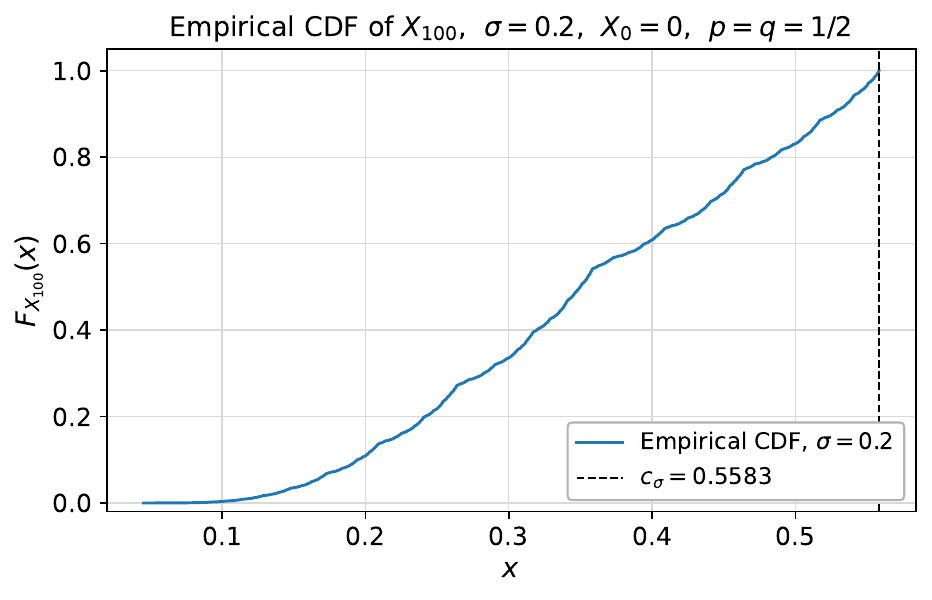}
\caption{Numerical illustration of the support geometry in the
overlapping regime \(\sigma=0.2<\sigma_c\). }
\label{fig:sigma-0p2}
\end{figure}

\begin{figure}[H]
\centering
\includegraphics[width=0.75\textwidth]{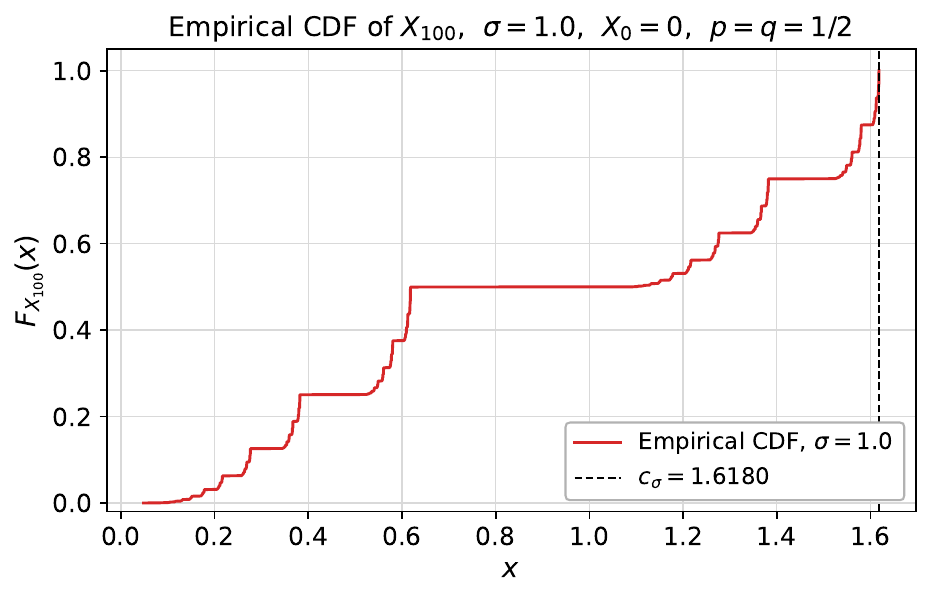}
\caption{Numerical illustration of the support geometry in the separated
regime \(\sigma=1>\sigma_c\).}
\label{fig:sigma-1p0}
\end{figure}
Figures~\ref{fig:sigma-0p2} and~\ref{fig:sigma-1p0} provide numerical
illustrations of the support geometry for \(\sigma=0.2\) and \(\sigma=1\).
Since \(X_0=0\) and the innovations are binary, the law of \(X_{100}\) is
discrete, and its empirical distribution function is necessarily a step
function regardless of the type of the limiting stationary distribution.
These plots are therefore intended as visual aids only.

\section{Conclusion}
The present paper studies a nonlinear autoregressive process with binary innovations, generated by a strictly increasing and strictly concave branch with a neutral fixed point. Despite the failure of uniform contractivity, a two-step synchronous coupling yields existence and uniqueness of the stationary probability measure, together with quantitative convergence in the 1-Wasserstein metric for initial distributions with finite first moments. Moreover, weak convergence to the unique stationary probability measure is established for arbitrary initial distributions. It is further shown that the noise amplitude governs a sharp geometric transition of the stationary support between a full interval and a Cantor set. In the separated regime, the support has zero Lebesgue measure and the stationary distribution is non-atomic, and hence singular continuous. These results illustrate how the interaction between neutral deterministic dynamics and discrete innovations can produce both quantitative stochastic stability and a pronounced change in the geometry and regularity of the stationary law.

\section*{Disclosure statement}

The author declares no competing interests.

\section*{Funding}

This research received no specific grant from any funding agency in the
public, commercial, or not-for-profit sectors.

\section*{Declaration of Generative AI Use}
During the preparation of this manuscript, the author used OpenAI ChatGPT
(GPT-5.6 Sol) for text editing and polishing, LaTeX-related assistance, and
debugging of the Python code used to generate the numerical illustrations.
The mathematical ideas, results, and proof strategies were developed by the
author, while ChatGPT assisted with the presentation and clarification of
selected arguments and calculations. All AI-assisted material was carefully
reviewed, verified, and, where necessary, corrected by the author, who takes
full responsibility for the content of the paper.

\bibliographystyle{amsplain}
\bibliography{biblio}

\end{document}